\documentclass[draft,reqno,11pt,oneside]{article} 

\usepackage{amssymb}
\usepackage{english}
\usepackage{amsthm}
\usepackage{amssymb}
\usepackage{longtable} 
\usepackage{amsmath}
\usepackage{ifthen} 
\usepackage{upref}
\usepackage{bbold} 
\usepackage{latexsym} 
\usepackage[all]{xy}
\usepackage{fancyhdr} 
\usepackage{makeidx} 
\newfont{\chen}{cmex10 at 9pt} \newfont{\chee}{cmex10 at 10pt} 
\newfont{\cheu}{cmex10 at 11pt} 
\makeindex 
\begin{document} 
\newtheorem{lemma}{Lemma}[section]
\newtheorem*{eldemo}{Sketch of proof}
\newtheorem{satz}[lemma]{Theorem}
\newtheorem{prop}[lemma]{Proposition}
\newtheorem{defi}[lemma]{Definition}
\newtheorem{bei}[lemma]{Example} 
\newtheorem{verm}[lemma]{Conjecture} 
\newtheorem{kor}[lemma]{Corollary} 
\renewcommand{\proofname}{Proof} 
\newtheorem{bez}[lemma]{Notation} 
\newtheorem{bem}[lemma]{Remark}
\newtheorem{fall}[lemma]{Case}
\newtheorem*{example}{Example} 
\newtheorem{defsatz}[lemma]{Definition and Theorem} 
%
\newtheorem{lemmaap}{Lemma}[section]
\newtheorem{satzap}[lemmaap]{Theorem}
\newtheorem{propap}[lemmaap]{Proposition}
\newtheorem{defiap}[lemmaap]{Definition}
\newtheorem{beiap}[lemmaap]{Example}  
\newtheorem{korap}[lemmaap]{Corollary}  
\newtheorem{bezap}[lemmaap]{Notation} 
\newtheorem{bemap}[lemmaap]{Remark}
\def \Q{{\mathbb{Q}}}
\def \Ar{{\mathcal{A}}}
\def \G{{\mathcal{G}}} 
\def \pu{{\mbox{.}}} 
\def \ps{{\Phi^s}} 
\def \LUC{{\mathrm{LUC}}} 
\def \RUC{{\mathrm{RUC}(\G)}}
\def \RSp{{\widetilde{R}}^{(p)}} 
\def \RSq{{\widetilde{R}}^{(q)}} 
\def \wwp{{\mathfrak{P}}} 
\def \phii{{\overline{\Phi}}} 
\def \Mun{{\mathcal{M}}}
\def \Ce{{\mathrm{C}}}
\def \Cede{{\mathbb{C}}} 
\def \lpp{{L_p}} 
\def \Te{{\mathcal{T}}}
\def \liw{{L_\infty \left( \G , w^{-1} \right)}} 
\def \sis{{\mathcal{S}}} 
\def \th{{\widetilde{\Theta}}} 
\def \mc{{\mathcal{MC}}} 
\def \ml{{\mathcal{MLUC}}} 
\def \sus{{\subseteq}} 
\def \mr{{\mathcal{MRUC}}} 
\def \N{{\mathbb{N}}}
\def \al{{\alpha}} 
\def \a{{\alpha}} 
\def \suppp{{\mathrm{supp}}} 
\def \cb{{\mathrm{cb}}} 
\def \A{{\mathcal{A}}}
\def \llu{{\mathrm{LUC}(\G)}} 
\def \B{{\mathcal{B}}}
\def \kapp{{\mathfrak{k}}} 
\def \ru{{{\RUC}^*}} 
\def \L{{\mathfrak{L}}}
\def \EL{{\mathcal{L}}}
\def \gip{{\Gamma_{l, p}}} 
\def \UC{{{\mathrm{UC}}(\G)}} 
\def \gup{{{\widetilde{\Gamma_{l, p}}}}} 
\def \gurp{{\Gamma_{r, p}}} 
\def \gurrp{{{\widetilde{\Gamma_{r, p}}}}} 
\def \M{{\rm{M}}}
\def \en{{\mathcal{N}}}
\def \gr{{\widetilde{\Gamma_r}}} 
\def \gre{{\Gamma_r}} 
\def \em{{\mathcal{M}}} 
\def \emli{{\mathcal{ML}_\infty}} 
\def \ccc{{\mathcal{C}}}
\def \Ha{{\mathfrak{H}}}
\def \l1{{L_1(\G)}} 
\def \li{{L_{\infty}(\G)}} 
\def \H{{\mathcal{H}}}
\def \Qu{{\mathfrak{Q}}}
\def \cbgl{~{\stackrel{\mathrm{cb}}{=}}~} 
\def \FS{{\mathcal{FS}}} 
\def \gammas{{\widetilde{\gamma}}} 
\def \C{{\mathcal{C}}}
\def \mr{{\mathcal{MRUC}}}
\def \Be{{\mathfrak{B}}}
\def \bar{{~|~}} 
\def \Ka{{\mathcal{K}}}
\def \rr{{\mathcal{R}}} 
\def \ceb{{\C\B(\B(L_2(\G)))}} 
\def \id{{{\mathrm{id}}}} 
\def \cep{{\C\B(\B(L_p(\G)))}} 
\def \ep{{\B(\B(L_p(\G)))}}
\def \Is{{\widetilde{I}}} 
\def \conv{{\mathrm{conv}}} 
\def \convs{{\widetilde{\mathrm{conv}}}} 
\def \convsp{{\widetilde{\mathrm{conv_p}}}} 
\def \es{{\B(\B(L_2(\G)))}}
\def \luw{{{\mathrm{LUC} \left( \G , w^{-1} \right)}^*}} 
\def \lucw{{\mathrm{LUC} \left( \G , w^{-1} \right)}} 
\def \c0w{{\Cee_0 \left( \G , w^{-1} \right)}} 
\def \mw{{\mathrm{M} (\G , w)}} 
\def \go{{\gip}} 
\def \deltas{{\widetilde{\delta}}} 
\def \tt{{\overline{\otimes}}} 
\def \ti{{\stackrel{\vee}{\otimes}}} 
\def \oo{{{\overline{\otimes}}}} 
\def \tp{{\widehat{\otimes}}} 
\def \otp{{\otimes}} 
\def \zet{{Z_t(\lu)}} 
\makeatletter 
\newcommand{\essu}
{\mathop{\operator@font ess\mbox{-}sup}} 
\newcommand{\essi}
{\mathop{\operator@font ess\mbox{-}inf}} 
\newcommand{\lisu}
{\mathop{\operator@font lim\,ess\mbox{-}sup}} 
\newcommand{\liin}
{\mathop{\operator@font lim\,ess\mbox{-}inf}} 
\makeatother 
\def \ga{{\Gamma_p}}
\def \gam{{\gamma_p}}
\def \gaw{{\widetilde{\ga}}}
\def \lu{{\LUC(\G)^*}}
\def \sa{{\overline{\Gamma_p(\M(\G))}^{w*}}}
\def \Ball{{\mathrm{Ball}}} 
\def \gurr{{\widetilde{\Gamma_r}}} 
\def \gur{{\Gamma_r}} 
\def \Ker{{\it{KERN}}}
\def \F{{\mathfrak{F}}}
\def \ad{{\mathrm{ad}}} 
\def \Ceee{{\mathbb{C}}}
\def \Cee{{\mathrm{C}}}
\def \imink{{\iota_{\mathrm{minK}}}} 
\def \imin{{\iota_{\mathrm{min}}}} 
\def \Bild{{\it{BILD}}}
\def \Sz{{\mathcal{S}}}
\def \The{{\mathcal{T}}}
\def \cebe{{\C\B(\B(\Sz_2))}}
\def \cebes{{\C\B(\B(\Sz_2(L_2(G))))}} 
\def \pr{{\Sz_2 \otimes_h \B(\Sz_2) \otimes_h \Sz_2}} 
\def \su{{\rm{sup}}_{t \in \G}}
\def \mi{{\Gamma_2(\M(\G))}}
\def \Cb{{\mathcal{C}}}
\newcommand{\ens}{{\mathcal{N}}(L_p(\G))}
\newcommand{\Tee}{{\mathcal{T}}(\mathcal{H})}
\newcommand{\schutz}[1]{#1} 
\def \ma{{\Gamma_2(\lu)}}
\def \Hi{{\mathcal{H}}}
\def \bs{{\B^{\sigma}(\B(\Ha))}}
\def \bsa{{\C\B(\B(\H))}}
\def \bsi{{\B^s(\B(\Ha))}}
\def \gu{{\widetilde{\Gamma_l}}}
\def \EF{{\widetilde{F}}} 
\def \l1w{{L_1(\G , w)}} 
\def \gi{{\Gamma_l}} 
\def \als{{\widetilde{\alpha}}} 
\def \betas{{\widetilde{\beta}}} 
\def \gir{{\Gamma_r}} 
\def \muss{{\widetilde{\mu}}} 
\def \ge{{\Gamma}}
\def \gem{{\Gamma(\M(\G))}}
\def \MM{{\mathbf{M}}} 
\def \un{{\ell_{\infty}^*(\G)}}
\def \R{{\mathcal{R}}}
\def \cebr{{\C\B_{\R(\G)}(\B(L_2(\G)))}}
\def \Chi{{\chi}} 
\def \cebre{{\C\B_{\R(\G)}(\B(\ell_2(\G)))}}
\def \cebrs{{\C\B_{\R(\G)^{'}}^{\sigma}(\B(L_2(\G)))}}
\def \cebra{{\C\B_{\R(\G)}^{\sigma}(\B(\ell_2(\G)))}}
\def \U{{\mathfrak{U}}}
\def \K{\Ka} 
\def \cebru{{\C\B_R(\B(L_2(\G)))}}
\def \lin{{\mathrm{lin}}} 
\def \cebri{{\C\B_{R}^{\sigma}(\B(\H))}}
\def \Tee{{\The(\H)}}
\def \gurrn{{\gurr^0}} 
\def \gin{{\gi^0}} 
\def \gun{{\gu^0}} 
\newcommand{\dogl}{
\setlength{\unitlength}{0.7ex}
\linethickness {0.1ex}
~\begin{picture} (2.90 , 2)
\put (0 , 0.5) {\circle* {0.3}}
\put (0 , 1.2) {\circle* {0.3}}
\put (0.7 , 0.5) {\line (1,0) {2.1}}
\put (0.7 , 1.2) {\line (1,0) {2.1}}
\end{picture}~}
%
%
\newlength{\checklength}
\newcommand{\fns}{\footnotesize}
\newcommand{\widecheck}[1]{
\settowidth{\checklength}{$#1$}
\linethickness {0.1ex}
\setlength{\unitlength}{\checklength}
\stackrel{
\ifthenelse{ \lengthtest{ \checklength > 6em } }
{
\begin{picture}(1,0.04)
\qbezier(0,0.04)(0.3,0.03)(0.5,0.003)
\qbezier(0.5,0.003)(0.7,0.03)(1,0.04)
\qbezier(0,0.04)(0.3,0.03)(0.5,-0)
\qbezier(0.5,-0)(0.7,0.03)(1,0.04)
\qbezier(0,0.04)(0.3,0.03)(0.5,-0.005)
\qbezier(0.5,-0.005)(0.7,0.03)(1,0.04)
\end{picture}
}
{
\ifthenelse{ \lengthtest{ \checklength > 4em } }
{
\begin{picture}(1,0.05)
\qbezier(0,0.05)(0.3,0.03)(0.5,-0.02)
\qbezier(0.5,-0.02)(0.7,0.03)(1,0.05)
\qbezier(0,0.05)(0.3,0.03)(0.5,-0.015)
\qbezier(0.5,-0.015)(0.7,0.03)(1,0.05)
\qbezier(0,0.05)(0.3,0.03)(0.5,-0.01)
\qbezier(0.5,-0.01)(0.7,0.03)(1,0.05)
\end{picture}
}
{
\ifthenelse{\lengthtest{\checklength > 3em}}
{
\begin{picture}(1,0.075)
\qbezier(0,0.075)(0.3,0.05)(0.5,-0.03)
\qbezier(0.5,-0.03)(0.7,0.05)(1,0.075)
\qbezier(0,0.075)(0.3,0.05)(0.5,-0.02)
\qbezier(0.5,-0.02)(0.7,0.05)(1,0.075)
\end{picture}
}
{
\ifthenelse{\lengthtest{\checklength > 2em}}
{
\begin{picture}(1,0.1)
\qbezier(0,0.1)(0.3,0.07)(0.5,-0.04)
\qbezier(0.5,-0.04)(0.7,0.07)(1,0.1)
\qbezier(0,0.1)(0.3,0.07)(0.5,-0.02)
\qbezier(0.5,-0.02)(0.7,0.07)(1,0.1)
\end{picture}
}
{
\ifthenelse{\lengthtest{\checklength > 1em}}
{
\begin{picture}(1,0.12)
\qbezier(0,0.12)(0.3,0.1)(0.5,-0.06)
\qbezier(0.5,-0.06)(0.7,0.1)(1,0.12)
\qbezier(0,0.12)(0.3,0.1)(0.5,-0.025)
\qbezier(0.5,-0.025)(0.7,0.1)(1,0.12)
\end{picture}
}
{
\begin{picture}(1,0.15)
\qbezier(0,0.15)(0.3,0.12)(0.5,-0.1)
\qbezier(0.5,-0.1)(0.7,0.12)(1,0.15)
\qbezier(0,0.15)(0.3,0.12)(0.5,-0.04)
\qbezier(0.5,-0.04)(0.7,0.12)(1,0.15)
\end{picture}
}
}
}
}
}
}
{#1}
}
\newcommand{\swidecheck}[1]{
\settowidth{\checklength}{$_{#1}$}
\linethickness {0.1ex}
\setlength{\unitlength}{0.95\checklength}
\stackrel{
\ifthenelse{ \lengthtest{ \checklength > 6em } }
{
\begin{picture}(1,0.04)
\qbezier(0,0.04)(0.3,0.03)(0.5,0.003)
\qbezier(0.5,0.003)(0.7,0.03)(1,0.04)
\qbezier(0,0.04)(0.3,0.03)(0.5,-0)
\qbezier(0.5,-0)(0.7,0.03)(1,0.04)
\qbezier(0,0.04)(0.3,0.03)(0.5,-0.005)
\qbezier(0.5,-0.005)(0.7,0.03)(1,0.04)
\end{picture}
}
{
\ifthenelse{ \lengthtest{ \checklength > 4em } }
{
\begin{picture}(1,0.05)
\qbezier(0,0.05)(0.3,0.03)(0.5,-0.02)
\qbezier(0.5,-0.02)(0.7,0.03)(1,0.05)
\qbezier(0,0.05)(0.3,0.03)(0.5,-0.015)
\qbezier(0.5,-0.015)(0.7,0.03)(1,0.05)
\qbezier(0,0.05)(0.3,0.03)(0.5,-0.01)
\qbezier(0.5,-0.01)(0.7,0.03)(1,0.05)
\end{picture}
}
{
\ifthenelse{\lengthtest{\checklength > 3em}}
{
\begin{picture}(1,0.075)
\qbezier(0,0.075)(0.3,0.05)(0.5,-0.03)
\qbezier(0.5,-0.03)(0.7,0.05)(1,0.075)
\qbezier(0,0.075)(0.3,0.05)(0.5,-0.02)
\qbezier(0.5,-0.02)(0.7,0.05)(1,0.075)
\end{picture}
}
{
\ifthenelse{\lengthtest{\checklength > 2em}}
{
\begin{picture}(1,0.1)
\qbezier(0,0.1)(0.3,0.07)(0.5,-0.04)
\qbezier(0.5,-0.04)(0.7,0.07)(1,0.1)
\qbezier(0,0.1)(0.3,0.07)(0.5,-0.02)
\qbezier(0.5,-0.02)(0.7,0.07)(1,0.1)
\end{picture}
}
{
\ifthenelse{\lengthtest{\checklength > 1em}}
{
\begin{picture}(1,0.12)
\qbezier(0,0.12)(0.3,0.12)(0.5,-0.06)
\qbezier(0.5,-0.06)(0.7,0.12)(1,0.12)
\qbezier(0,0.12)(0.3,0.12)(0.5,-0.025)
\qbezier(0.5,-0.025)(0.7,0.12)(1,0.12)
\end{picture}
}
{
\begin{picture}(1,0.15)
\qbezier(0,0.15)(0.3,0.12)(0.5,-0.1)
\qbezier(0.5,-0.1)(0.7,0.12)(1,0.15)
\qbezier(0,0.15)(0.3,0.12)(0.5,-0.04)
\qbezier(0.5,-0.04)(0.7,0.12)(1,0.15)
\end{picture}
}
}
}
}
}
}
{#1}
} 
\def \des{{\sqrt{x_j} \otimes M_{T_{\widecheck{b}}} \otimes \sqrt{x_j}}} 
\def \res{{(e_j)^{\frac{1}{p}} \otimes T_{\widecheck{b}} \otimes (e_j)^{\frac{1}{q}}}} 
\def \rem{{(e_j)^{\frac{1}{2}} \otimes M_{\widecheck{b}} \otimes (e_j)^{\frac{1}{2}}}} 
%
%
\title{On the Topological Centre 
Problem \\ 
for Weighted Convolution Algebras} 
\author{Matthias Neufang}
\date{} 
\maketitle 
\begin{abstract} 
\noindent 
Let $\G$ be a locally compact non-compact group. 
We\footnote{ 
\par 
The author is currently a PIMS Postdoctoral Fellow at the 
University of Alberta, Edmonton, where this work was accomplished. The support of PIMS is gratefully acknowledged.} 
show that under a very mild assumption on the weight function $w$, the weighted group algebra $\l1w$ is strongly 
Arens irregular in the sense of \cite{dll}. To this end, we first derive a general factorization theorem for bounded 
families in the $\liw^*$-module $\liw$. 
\end{abstract}

\section{Introduction} 
Let $\G$ be a locally compact group, 
and let $w: \G \longrightarrow ( 0 , \infty )$ be a weight function, i.e., a positive continuous 
function on $\G$ such that $w(st) \leq w(s) w(t)$ for all $s,t \in \G$; for convenience we shall assume that $w(e)=1$, 
where $e$ is the neutral element of $\G$. We 
will consider the following spaces, normed in 
such a way that multiplication resp.\ division by the weight becomes an isometry between the unweighted and the 
corresponding weighted space (whose norm we will denote by $\| \cdot \|_w$): 
\begin{eqnarray*}
\l1w &=& \{ f \mid w f \in L_1(\G) \} \\ 
\liw &=&  \{ f \mid w^{-1} f \in L_{\infty}(\G) \} \\ 
\lucw &=& \{ f \mid w^{-1} f \in \LUC (\G) \} \\ 
\c0w &=& \{ f \mid w^{-1} f \in \Cee_0 (\G) \} \\ 
\mw &=& \{ \mu \mid w \mu \in \M(\G) \} . 
\end{eqnarray*} 
Then we have $\liw = \l1w^*$ and $\mw = \c0w^*$. 
For every $y \in \G$, we define $\deltas_{y} := w(y)^{-1} ~\delta_y$, 
which is an element of norm one in $\mw$. 
\par 
Our aim is to show that for all locally compact non-compact groups, the weighted group algebra $\l1w$ is 
strongly Arens irregular in the sense of Dales--Lamb--Lau (see \cite{dll}), provided 
the weight satisfies some 
very mild boundedness condition. Here, strong Arens irregularity means that 
the topological centre of the bidual algebra $( \l1w^{**}, \odot )$, equipped with the first Arens product, 
precisely equals the algebra $\l1w$ itself, i.e., it is extremally small. This is a generalization of the main result, Thm.\ 1, 
of \cite{lalo}, where the corresponding assertion is proved for the (unweighted) group algebra $L_1(\G)$, to the 
weighted situation. 
Although covering a by far more general case, our proof 
is not of higher complexity, if not even simpler, than the one given in \cite{lalo}. 
\par 
In the following, we shall always regard 
$\l1w^{**}$ as endowed with the first Arens multiplication. Let us briefly recall the three step construction of the latter, 
arising from the convolution product (denoted by ``$*$'') in $\l1w$ via various 
module actions. -- For 
$m , n \in \l1w^{**}$, $h \in \l1w^{*}$ and $f , g \in \l1w$ one defines: 
\begin{center} 
\begin{tabular}{ccc} 
$\langle h \odot f, g \rangle$ & $:=$ & $\langle h, f * g \rangle$ \\ 
$\langle n \odot h, f \rangle$ & $:=$ & $\langle n, h \odot f \rangle$ \\ 
$\langle m \odot n, h \rangle$ & $:=$ & $\langle m, n \odot h \rangle$. 
\end{tabular} 
\end{center} 
A fairly comprehensive ex\-position of the 
basic theory of Arens products is given in \cite{pal}, \S 1.4. As for topological centres, an excellent source is 
\cite{laul}. We shall only need the definition of the latter, which we briefly recall here: 
$$Z_t \left( \l1w^{**} \right) := \{ m \in \l1w^{**} \mid n \mapsto m \odot n ~\mathrm{is} ~w^*-w^*-\mathrm{continuous ~on} ~\l1w^{**} \} .$$ 
We will use the fact (cf.\ \cite{gro}, Prop.\ 1.3) that, 
with the natural module operation stemming from the construction of the first Arens 
product on $\l1w^{**}$, the equality $\liw \odot \l1w = \lucw$ holds. 
Hence, 
a natural 
module operation 
of $\luw$ on $\liw$ is given by 
$$\langle m \diamond h, g \rangle = \langle m, h \odot g \rangle ,$$ 
where $m \in \luw$, $h \in \liw$, $g \in \l1w$. 
It is readily verified that we have $m \diamond h = \widetilde{m} \odot h$, 
where $\widetilde{m}$ is an arbitrary Hahn-Banach extension of $m$ to $\liw^*$. 
\par 
In the sequel, we shall denote by $\kapp (\G)$ the \textit{compact covering number} of the group $\G$, i.e., 
the least cardinality of a compact covering of $\G$. 
For the sake of brevity, we further introduce the following terminology (the first part of the definition also 
appears in \cite{dll}). 
\begin{defi} \label{dia} 
(i) A subset $S$ of $\G$ will be called $\mathrm{dispersed}$ if $S$ is not contained in any union 
of a family of compact subsets of $\G$, the family having cardinality strictly less than $\kapp (\G)$. 
\par 
(ii) For a subset $S \subseteq \G$, we say that the weight $w$ is $\mathrm{diagonally ~bounded ~on}$ $S$ if we have: 
$$\sup_{s \in S} ~w(s) w \left( s^{-1} \right) < \infty .$$ 
\end{defi}
Now we can formulate the main result of the present note; we remark that it has very recently also been obtained independently 
by Dales--Lamb--Lau in \cite{dll}, though with a different proof. In particular, our factorization result, Theorem \ref{fakt}, 
does not appear in \cite{dll}. 
\begin{satz} \label{hau} 
Let $\G$ be a locally compact non-compact group with compact covering number $\kapp(\G)$. 
Suppose that there is a dispersed set $S \subseteq \G$ on which the weight function $w$ is diagonally bounded. 
Then $\l1w$ is strongly Arens irregular. 
\end{satz} 
We wish to stress the following important points 
regarding our approach: 
\begin{itemize}
\item 
We prove a 
(formal) sharpening of the interesting inclusion 
contained in Theorem \ref{hau}. Namely, we will show that for an element 
$m \in \liw^*$ in order to belong to $\l1w$, it suffices 
that left multiplication by $m$ be $w^*$-$w^*$-continuous on the $w^*$-closure of the set of all 
Hahn-Banach extensions of functionals in 
${\overline{\widetilde{\delta}_\G}}^{w^*} \subseteq \Ball(\luw)$ to $\liw^*$. Instead, the definition of the topological centre 
demands $w^*$-continuity 
on all of $\liw^*$. 
\item 
The proof is \textit{direct} and follows, once the necessary prerequisites are established (section 2), 
a purely Banach algebraic pattern (section 3). Our global procedure is similar to the one presented in 
\cite{ich1a} and may thus be considered as having the same merits as the latter. 
\end{itemize} 

\section{A factorization theorem for families of functions in $\liw$} 
In the proof of our main result 
we will use the following two results, which are of independent interest. 
\begin{prop} \label{kmazur}
For an arbitrary locally compact group $\G$, the space $\l1w$ enjoys Mazur's property of 
level $\kapp(\G) \cdot \aleph_0$. 
-- This means that 
a functional $m \in \l1w^{**}$ actually belongs to $\l1w$ if it carries 
bounded $w^*$-converging nets of cardinality at most $\kapp(\G) \cdot \aleph_0$ into converging nets. 
\end{prop}
\begin{proof}
This follows from Thm.\ 4.4 in \cite{ich1}, which states the above property for the space $L_1(\G)$, 
and the fact that the latter is stable under isomorphism (cf.\ \cite{ich1}, Remark 4.3). 
\end{proof} 
\noindent 
Next we present 
our crucial tool, which is a general factorization theorem for 
bounded families in $\liw$; it provides a generalization of Thm.\ 2.2 in \cite{ich1a} to the weighted situation. 
\begin{satz} \label{fakt}
Let $\G$ be a locally compact non-compact group with compact covering number $\kapp(\G)$. 
Suppose that there is a dispersed set $S \subseteq \G$ on which the weight function $w$ is diagonally bounded. 
Then there exists 
a family $(\psi_\al)_{\al \in I}$, $|I|= \kapp(\G)$, 
of functionals in ${\overline{\widetilde{\delta}_\G}}^{w^*} \subseteq \Ball \left( \luw \right)$ such that 
whenever 
$(h_\al)_{\al \in I} \subseteq \liw$ is a bounded family of functions, there exists 
a function $h \in \liw$ such that the factorization formula 
$$h_\al = \psi_\al \diamond h$$
holds for all $\al \in I$. 
\end{satz}
\begin{proof} 
For $y \in \G$, we denote by $r_y$ the operator of right translation, i.e., $(r_y f) (x)=f(xy)$ whenever $f$ is a function on $\G$ 
and $x \in \G$. 
\par 
There is a covering of $\G$ by open sets whose closure is compact, 
of minimal cardinality, i.e., of cardinality $\kapp(\G)$, and closed under finite unions; 
we denote the corresponding family of compacta by $(K_\al)_{\al \in I}$. Set $\Is := I \times I$. 
For $\als = (\al , i) \in \Is$, put $K_\als = K_{(\al , i)} := K_\al$. Then 
$(K_\als)_{\als \in \Is}$ is a covering of $\G$ having the same properties than the original one. 
Since the set $S$ is dispersed, by the same reasoning as in Lemma 3 of \cite{lalo}, we see that 
there exists a family $(y_\als)_{\als \in \Is} \subseteq S$ such that 
\begin{eqnarray} \label{0} 
K_\als y_\als^{-1} \cap K_\betas y_\betas^{-1} = \varnothing \quad \forall \als , \betas \in \Is , ~\als \not= \betas . 
\end{eqnarray} 
Set $S' := \left\{ y_\als \mid \als \in \Is \right\}$. We define, for $(\al , i) , (\beta , j) \in \Is$: 
\begin{eqnarray} \label{1} 
(\al , i) \preceq (\beta , j) :\Longleftrightarrow K_{(\al , i)} \subseteq K_{(\beta , j)} 
\Longleftrightarrow K_\al \subseteq K_\beta \Longleftrightarrow : \al \preceq' \beta . 
\end{eqnarray} 
Let $\F$ be an ultrafilter on $I$ which dominates the order filter. Define, for $j \in I$, 
$${\psi_{j}}' := w^*-\lim_{\beta \rightarrow \F} \deltas_{y^{-1}_{(\beta , j)}} \in {\overline{{\deltas_\G}}}^{w^*} 
\subseteq \Ball \left( \luw \right) ,$$ 
and let $\psi_j$ be arbitrary Hahn-Banach extensions of ${\psi_{j}}'$ to $\liw^*$. 
\par 
Since $w$ is diagonally bounded on $S'$, we have: 
$$\sup_{s \in S'} ~\left\| w \left( s^{-1} \right) \delta_s \right\|_w = \sup_{s \in S'} ~w(s) w \left( s^{-1} \right) < \infty .$$ 
Thus, the family of functions 
$$H_{(\al ,i)} := \left( w \left( y^{-1}_{(\al , i)} \right) ~\delta_{y_{(\al , i)}} \right) \diamond 
\left( \Chi_{K_{(\al ,i)}} ~h_i \right) = w \left( y^{-1}_{(\al , i)} \right) ~r_{y_{(\al , i)}} ~\left( 
\Chi_{K_{(\al ,i)}} ~h_i \right)$$ 
is bounded in $\liw$, whence $\left( w^{-1} H_{(\al ,i)} \right)$ is a bounded family in $\li$. 
By (\ref{0}), the projections $r_{y_{(\al , i)}} ~\Chi_{K_{(\al ,i)}} = \Chi_{K_{(\al ,i)} y_{(\al , i)}^{-1}}$ 
are pairwise orthogonal, so that 
$$H := \sum_{\al \in I} \sum_{i \in I} w^{-1} H_{(\al ,i)} \quad (w^*-\mathrm{limits})$$ 
defines a function in $L_\infty (\G)$. Hence, we have 
$$h := \sum_{\al \in I} \sum_{i \in I} H_{(\al ,i)} \in \liw .$$ 
Using (\ref{0}), 
we obtain for all $(\al , i)$, $(\beta , j)$, $(\gamma , k) \in \Is$, where $(\gamma , k) \preceq (\beta , j)$: 
\begin{eqnarray*} 
\Chi_{K_{(\gamma , k)}} ~r_{y^{-1}_{(\beta , j)}} r_{y_{(\al , i)}} ~\left( \Chi_{K_{(\al ,i)}} ~h_i \right) &=& 
\Chi_{K_{(\gamma , k)}} ~\Chi_{K_{(\beta , j)}} ~r_{y^{-1}_{(\beta , j)}} r_{y_{(\al , i)}} ~\left( \Chi_{K_{(\al ,i)}} ~h_i \right) \\ 
&=& \Chi_{K_{(\gamma , k)}}~ \left[ r_{y^{-1}_{(\beta , j)}} \left(  \left(  r_{y_{(\beta , j)}} ~\Chi_{K_{(\beta , j)}}  \right)  
~r_{y_{(\al , i)}} ~\left( \Chi_{K_{(\al ,i)}} ~h_i  \right)  \right) \right] \\ 
&=& \delta_{(\al , i),(\beta , j)} ~\Chi_{K_{(\gamma , k)}} ~h_j . 
\end{eqnarray*} 
Taking into account (\ref{1}), we deduce that for all 
$j \in I$ and $(\gamma , k) \in \Is$: 
\begin{eqnarray*} 
\Chi_{K_{(\gamma , k)}} ~\left( \psi_j \diamond h \right) &=& 
w^*-\lim_{\beta \rightarrow \F} ~\sum_{\al \in I} \sum_{i \in I} 
w \left( y^{-1}_{(\al , i)} \right) w \left( y^{-1}_{(\beta , j)} \right)^{-1} 
~\underbrace{\Chi_{K_{(\gamma , k)}} ~r_{y^{-1}_{(\beta , j)}} r_{y_{(\al , i)}} ~\left( \Chi_{K_{(\al ,i)}} ~h_i \right)}_{= 
~\delta_{(\al , i),(\beta , j)} ~\Chi_{K_{(\gamma , k)}} ~h_j} \\ 
&=& \Chi_{K_{(\gamma , k)}} ~h_j , 
\end{eqnarray*} 
whence the desired factorization formula follows. 
\end{proof} 

\section{Strong Arens irregularity of $\l1w$} 
We now come to 
the proof of Theorem \ref{hau}. -- To establish the 
non-trivial inclusion, let $m \in Z_t \left( \l1w^{**} \right)$. The group $\G$ being non-compact, we infer from 
Proposition \ref{kmazur} that $\l1w$ has Mazur's property of level $\kapp(\G)$. So in order to prove that $m \in \l1w$, let 
$(h_\al)_{\al \in I} \subseteq \liw$ be a bounded net converging $w^*$ to $0$, where $|I|=\kapp(\G)$. 
Thanks to Theorem \ref{fakt}, we have the factorization 
$$h_\al = \psi_\al \diamond h = 
\widetilde{\psi_\al} \odot h \quad (\al \in I)$$ 
with $\psi_\al \in {\overline{\widetilde{\delta}_\G}}^{w^*} \subseteq \Ball(\luw)$ and $h \in \liw$. -- Here, 
$\widetilde{\psi_\al}$ denotes some arbitrarily chosen Hahn-Banach extension of $\psi_\al$ to $\liw^*$. 
We have to show that $a_\al := \langle m, h_\al \rangle \stackrel{\al}{\longrightarrow} 0$. Due to the 
boundedness of $(h_\al)_\al$, it suffices to prove that every convergent subnet of $(a_\al)_\al$ tends to $0$. Let 
$(\langle m, h_{\al_\beta} \rangle)_\beta$ be such a convergent subnet. Furthermore, let 
$$E := w^*-\lim_\gamma \widetilde{\psi_{\al_{\beta_\gamma}}} \in \Ball \left( \liw^* \right)$$ 
be a $w^*$-cluster point of the net ${\left( \widetilde{\psi_{\al_{\beta}}} \right)}_\beta \subseteq \Ball \left( \liw^* \right)$. 
\par 
We first note that 
$E \odot h=0$, since for arbitrary $g \in \l1w$ we obtain: 
\begin{eqnarray*} 
\langle E \odot h, g \rangle &=& \langle E, h \odot g \rangle = 
\lim_\gamma \langle \psi_{\al_{\beta_\gamma}}, h \odot g \rangle \\ 
&=& \lim_\gamma \langle \psi_{\al_{\beta_\gamma}} \diamond h, g \rangle 
= \lim_\gamma \langle h_{\al_{\beta_\gamma}}, g \rangle \\ 
&=& 0. 
\end{eqnarray*}
Now we conclude, using the fact that $m \in Z_t \left( \l1w^{**} \right)$: 
\begin{eqnarray*} 
\lim_\beta \langle m, h_{\al_\beta} \rangle &=& 
\lim_\gamma \langle m, h_{\al_{\beta_{\gamma}}} \rangle = \lim_\gamma \langle m, \widetilde{\psi_{\al_{\beta_\gamma}}} \odot h \rangle \\ 
&=& \lim_\gamma \langle m \odot \widetilde{\psi_{\al_{\beta_\gamma}}}, h \rangle 
= \langle m \odot E , h \rangle \\ 
&=& \langle m , E \odot h \rangle = 0 , 
\end{eqnarray*} 
which yields the desired convergence.

\vspace{0.3cm} 
\noindent 
{\sc{Author's address:}} 
\\ 
{\textit{
Department of Mathematical Sciences\\ 
University of Alberta\\ 
Edmonton, Alberta\\ 
Canada T6G 2G1\\ 
E-mail:}} 
mneufang@math.ualberta.ca 
\end{document}